\documentclass[12pt]{amsart}

\usepackage{amssymb,amsmath,amsthm,enumitem,hyperref}
\setlist[enumerate,1]{label=\textup{(\alph*)},ref=(\alph*)}

\newcommand{\NN}{\mathbb{N}}
\newcommand{\PP}{\mathbb{P}}

\theoremstyle{plain}
\newtheorem*{conjecture2}{Demailly's Conjecture}
\newtheorem{theorem}{Theorem}
\newtheorem{corollary}[theorem]{Corollary}

\begin{document}

\mbox{}
\vspace{-1.1ex}
\title[Demailly's conjecture on Waldschmidt constants for many points]{Demailly's conjecture on Waldschmidt constants for sufficiently many very general points in $\PP^n$}
\author{Yu-Lin Chang}
\address{Department of Mathematics \\
National Tsing Hua University \\
Taiwan}
\email{\texttt{s103021505@m103.nthu.edu.tw}}
\author{Shin-Yao Jow}
\address{Department of Mathematics \\
National Tsing Hua University \\
Taiwan}
\email{\texttt{syjow@math.nthu.edu.tw}}
\date{}

\begin{abstract}
Let $Z$ be a finite set of $s$ points in the projective space $\PP^n$ over an algebraically closed field $F$. For each positive integer $m$, let $\alpha(mZ)$ denote the smallest degree of nonzero homogeneous polynomials in $F[x_0,\ldots,x_n]$ that vanish to order at least $m$
at every point of $Z$. The Waldschmidt constant $\widehat{\alpha}(Z)$ of $Z$ is defined by the limit \[
  \widehat{\alpha}(Z)=\lim_{m \to \infty}\frac{\alpha(mZ)}{m}. \]
Demailly conjectured that
\[
\widehat{\alpha}(Z)\geq\frac{\alpha(mZ)+n-1}{m+n-1}.
\]
Malara, Szemberg, and Szpond~\cite{gM17} established Demailly's conjecture when $Z$ is very general and \[
 \lfloor\sqrt[n]{s}\rfloor-2\geq m-1. \]
Here we improve their result and show that Demailly's conjecture holds if $Z$ is very general and \[
\lfloor\sqrt[n]{s}\rfloor-2\ge \frac{2\varepsilon}{n-1}(m-1), \]
where $0\le \varepsilon<1$ is the fractional part of $\sqrt[n]{s}$. In particular, for $s$ very general points where $\sqrt[n]{s}\in\NN$ (namely $\varepsilon=0$), Demailly's conjecture holds for all $m\in\NN$. We also show that Demailly's conjecture holds if $Z$ is very general and \[
 s\ge\max\{n+7,2^n\}, \]
assuming the Nagata-Iarrobino conjecture $\widehat{\alpha}(Z)\ge\sqrt[n]{s}$.
\end{abstract}

\keywords{Waldschmidt constant, Demailly's conjecture, Chudnovsky's conjecture, Nagata-Iarrobino conjecture}
\maketitle

\section{Introduction} 
Let $Z=\{p_1,\ldots,p_s\}$ be a set of $s$ points in the projective space $\PP^n$ over an algebraically closed field $F$, and let $I\subseteq F[x_0,\ldots,x_n]$ be the homogeneous ideal of $Z$. For any homogeneous ideal $J\subseteq F[x_0,\ldots,x_n]$, we write \[
 \alpha(J)=\min\{\deg f\mid f\text{ is a nonzero homogeneous polynomial in }J\}. \]
Given a positive integer $m$, an interesting question is how large the degree of a hypersurface has to be in order for it to pass through each point of $Z$ with multiplicity at least $m$. In more algebraic terms, this is asking about $\alpha(I^{(m)})$, where \[
 I^{(m)}=\langle f\in F[x_0,\ldots,x_n]\mid f\text{ is homogeneous and vanishes to order $\ge m$ at each }p_i \rangle \]
is the $m$-th symbolic power of $I$.

An asymptotic invariant closely related to this question is the \emph{Waldschmidt constant} $\widehat{\alpha}(I)$ of $I$, defined by \[
\widehat{\alpha}(I)=\inf_{m\ge 1}\frac{\alpha(I^{(m)})}{m}. \]
It is known \cite[Lemma~2.3.1]{cB01} that, in fact, \[
 \widehat{\alpha}(I) = \lim_{m \to \infty}\frac{\alpha(I^{(m)})}{m}.\]
Waldschmidt constants have been studied in various branches of mathematics, such as complex analysis \cite{hS75,jM80}, commutative algebra \cite{bH01}, and algebraic geometry \cite{hE83}. By definition, $\widehat{\alpha}(I)\le \alpha(I^{(m)})/m$ for all $m\in\NN$. An interesting conjecture of Demailly predicts an inequality in the reverse direction.

\begin{conjecture2}[\cite{jP82}]
Let $I$ be the homogeneous ideal of any finite set of $s$ points in $\PP^n$. Then \[
\widehat{\alpha}(I)\geq\frac{\alpha(I^{(m)})+n-1}{m+n-1} \]
for all $m\in\NN$.
\end{conjecture2}

Demailly's conjecture is known for $n=2$ over an algebraically closed field of characteristic $0$ \cite{hE83}. Demailly's conjecture for $m=1$ (also called the Chudnovsky's conjecture) is known for very general points over an algebraically closed field of characteristic $0$ \cite{lF18}, and for $s\geq 2^n$ very general points over an arbitrary algebraically closed field \cite{mD17}.

Recently, Malara, Szemberg, and Szpond~\cite{gM17} established Demailly's conjecture for $s\ge (m+1)^n$ very general points in $\PP^n$.\footnote{There is a preprint that improves this to $s\ge m^n$: see \cite[Theorem~4.8]{mD18}.} For later comparison with our Corollary~\ref{c:no Nagata}, note that the condition $s\ge (m+1)^n$ is equivalent to $\lfloor\sqrt[n]{s}\rfloor\ge m+1$, or \[
 \lfloor\sqrt[n]{s}\rfloor-2\ge m-1. \]

We now state our main theorem in this paper.

\begin{theorem}\label{t:main}
Let $I$ be the homogeneous ideal of any finite set of $s$ points in $\PP^n$.
\begin{enumerate}
  \item\label{unconditional} For all $m\in\mathbb{N}$, \[
    \frac{\alpha(I^{(m)})+n-1}{m+n-1}\le \sqrt[n]{s}. \]
  \item\label{conditional} If $\lfloor\sqrt[n]{s}\rfloor-2\ge \frac{2\varepsilon}{n-1}(m-1)$,
   where $0\le \varepsilon<1$ is the fractional part of $\sqrt[n]{s}$, then \[
    \frac{\alpha(I^{(m)})+n-1}{m+n-1}\le \lfloor\sqrt[n]{s}\rfloor. \]
\end{enumerate}
\end{theorem}

Since $\widehat{\alpha}(I)\ge\lfloor\sqrt[n]{s}\rfloor$ if $I$ is the homogeneous ideal of $s$ very general points in $\PP^n$ \cite{eL05,mD17}, Theorem~\ref{t:main}~\ref{conditional} implies the following Corollary~\ref{c:no Nagata}, which is an improvement of the aforementioned result in \cite{gM17} by Malara, Szemberg, and Szpond.

\begin{corollary}\label{c:no Nagata}
Demailly's conjecture holds for $s$ very general points in $\PP^n$ as long as \[
\lfloor\sqrt[n]{s}\rfloor-2\ge \frac{2\varepsilon}{n-1}(m-1), \]
where $0\le \varepsilon<1$ is the fractional part of $\sqrt[n]{s}$.
In particular, for $s$ very general points where $\sqrt[n]{s}\in\NN$ (namely $\varepsilon=0$), Demailly's conjecture holds for all $m\in\NN$.
\end{corollary}

It was conjectured by Iarrobino \cite{tI97} (and by Nagata \cite{mN59} for $n=2$) that $\widehat{\alpha}(I)\ge\sqrt[n]{s}$ if $I$ is the homogeneous ideal of $s\ge\max\{n+7,2^n\}$ very general points in $\PP^n$. Hence Theorem~\ref{t:main}~\ref{unconditional} implies the following

\begin{corollary}\label{c:Nagata}
Demailly's conjecture holds for $s\ge\max\{n+7,2^n\}$ very general points in $\PP^n$, provided that the Nagata-Iarrobino conjecture holds.
\end{corollary}

\section{Proof of Theorem~\ref{t:main}}
\ref{unconditional} Let $\delta=\lfloor \sqrt[n]{s}(m+n-1)-n+1\rfloor.$
Since $0\leq(\sqrt[n]{s}(m+n-1)-n+1)-\delta<1$,
\begin{align*}
\delta &>\sqrt[n]{s}(m+n-1)-n+1-1\\
&= \sqrt[n]{s}\cdot m + (\sqrt[n]{s}-1)(n-1)-1 \\
& \geq \sqrt[n]{s}\cdot m + (\sqrt[n]{s}-1)(n'-1)-1 &
&\quad\text{for $n' = 1,\ldots,n$}\\
& = \sqrt[n]{s}(m+n'-1)-n'+1-1  \\
& = \sqrt[n]{s}(m+n'-1)-n'.
\end{align*}
Therefore, we have
\begin{align*}
\delta +n' &> \sqrt[n]{s}(m+n'-1)
&\quad\text{for $n' = 1,\ldots,n$}
\end{align*}
and then
\[
 (\delta +n)(\delta +n-1)\ldots(\delta +1) > s(m+n-1)(m+n-2)\ldots m,
\]
which implies
\[
\binom{\delta+n}{n}>\binom{m+n-1}{n}\cdot s.
\]
It thus follows from dimension count that there exists a homogeneous polynomial of degree $\delta$ in $F[x_0,\ldots,x_n]$ vanishing at any given $s$ points in $\PP^n$ to order at least $m$. Hence \[
\alpha(I^{(m)})\le\delta=\lfloor\sqrt[n]{s}(m+n-1)\rfloor-n+1.\]
Thus \[
\frac{\alpha(I^{(m)})+n-1}{m+n-1}\leq \frac{\lfloor\sqrt[n]{s}(m+n-1)\rfloor}{m+n-1}\leq\sqrt[n]{s}. \]

\ref{conditional}
Let $k=\lfloor\sqrt[n]{s}\rfloor\in\mathbb{N}$. We want to show that if
\begin{equation}\label{E:assumption}
  k-2\ge \frac{2\varepsilon}{n-1}(m-1),
\end{equation}
then $\alpha(I^{(m)})\leq k(m+n-1)-n+1$, that is, $I^{(m)}$ contains a homogeneous polynomial of degree $k(m+n-1)-n+1$. By dimension count, it suffices to show that
\begin{equation}\label{E:sufficient condition}
\binom{k(m+n-1)+1}{n}>\binom{m+n-1}{n}\cdot s.
\end{equation}

Since $\sqrt[n]{s}=k+\varepsilon$, the inequality \eqref{E:sufficient condition} can be written as
\[
\binom{k(m+n-1)+1}{n}>\binom{m+n-1}{n}\cdot (k+\varepsilon)^n,
\]
which is equivalent to
\[
\bigl(k(m+n-1)+1\bigr)\cdots\bigl(k(m+n-1)+1-(n-1)\bigr)> (m+n-1)\cdots m\cdot (k+\varepsilon)^n.
\]
So to establish \eqref{E:sufficient condition}, it is sufficient to show that
\begin{multline}\label{E:i-th inequ}
\bigl(k(m+n-1)+1-i\bigr)\bigl(k(m+n-1)+1-(n-1)+i\bigr)\\
>(k+\varepsilon)^2(m+n-1-i)(m+i)
\end{multline}
for each $i=0,1,\ldots,\lfloor\frac{n-1}{2}\rfloor$.

Write the left-hand side minus the right-hand side of \eqref{E:i-th inequ} as a polynomial in $m$:\[
 \text{LHS\eqref{E:i-th inequ}}-\text{RHS\eqref{E:i-th inequ}}=Am^2+Bm+C_i, \]
where
\[
A=k^2-(k+\varepsilon)^2=-2k\varepsilon-\varepsilon ^2\le 0,
\]
\[
B=(n-1)(k^2-(1+2\varepsilon)k-\varepsilon ^2)+2k,
\]
and
\[
C_i=\bigl((k+\varepsilon)^2-1\bigr)i^2 - (n-1)\bigl((k+\varepsilon)^2-1\bigr)i + \bigl(k(n-1)+1\bigr)\bigl((k-1)(n-1)+1\bigr).
\]

Viewed as a quadratic function of $i$, $C_i$ attains its minimum at $i=\frac{n-1}{2}$, and the minimum value is \[
 C=\left(\frac{2k-1}{2}(n-1)+1\right)^2-\left(\frac{k+\varepsilon}{2}(n-1)\right)^2>0.
\]
Here $C>0$ because $k\ge 2$ by assumption \eqref{E:assumption}. Since $C_i\ge C$ for all $i$, to establish \eqref{E:i-th inequ}, it suffices to show that $Am^2+Bm+C>0$.

If $\varepsilon=0$, then $A=0$, $B>0$, and $C>0$, so $Am^2+Bm+C>0$ for all $m\in\NN$. If $\varepsilon>0$, then the assumption \eqref{E:assumption} can be rewritten as \[
 m\le \frac{n-1}{2\varepsilon}(k-2)+1. \]
To show that $Am^2+Bm+C$ is positive under this assumption, it suffices to show that $Am^2+Bm+C>0$ at $m=\frac{n-1}{2\varepsilon}(k-2)+1$, because $A<0$ and $C>0$.

Set $m=\frac{n-1}{2\varepsilon}(k-2)+1$, and denote $N=n-1$. Then \[
 Am^2+Bm+C=f(k,\varepsilon)N^2+g(k,\varepsilon)N+2k(1-\varepsilon)+(1-\varepsilon^2),\]
where \[
 f(k,\varepsilon)=\left(\frac{1}{2\varepsilon}-\frac{1}{2}\right)k^2+\left(-\frac{1}{\varepsilon}+2-\varepsilon\right)k-\frac{1}{4}(3-\varepsilon)(1-\varepsilon)\]
and \[
 g(k,\varepsilon)=\left(\frac{1}{\varepsilon}-1\right)k^2+\left(-\frac{2}{\varepsilon}+5-3\varepsilon\right)k-(1-\varepsilon)^2.\]

Note that $f(k,\varepsilon)$ and $g(k,\varepsilon)$ both tend to $0$ as $\varepsilon\to 1$, and their partial derivatives with respect to $\varepsilon$ are negative on $0<\varepsilon<1$, $k\geq2$, which implies that they are both positive on $0<\varepsilon<1$ and $k\geq2$: \[
\frac{\partial f}{\partial\varepsilon}=\frac{-1}{2\varepsilon^2}k^2+\left(\frac{1}{\varepsilon^2}-1\right)k-\frac{1}{4}(-4+2\varepsilon)\leq\frac{-1}{\varepsilon^2}k+\frac{1}{\varepsilon^2}k-k+1-\frac{\varepsilon}{2}<0,
\]
\[
\frac{\partial g}{\partial\varepsilon}=\frac{-1}{\varepsilon^2}k^2+\left(\frac{2}{\varepsilon^2}-3\right)k+2(1-\varepsilon)\leq\frac{-2}{\varepsilon^2}k+\frac{2}{\varepsilon^2}k-3k+2(1-\varepsilon)<0.
\]
Therefore $Am^2+Bm+C>0$ at $m=\frac{n-1}{2\varepsilon}(k-2)+1$, and the proof is complete.

\addcontentsline{toc}{section}{Reference}

\end{document}